\documentstyle[12pt]{article}
\input{amssym}
\newcommand{\g}{\goth g}
\begin{document}
\title{The special linear representations \\ of  compact Lie groups}
\author{M. Nadjafikhah\thanks{e-mail: m\_nadjafikhah@iust.ac.ir
Department of Mathematics, Iran University of Science and
Technology, Narmak, Tehran, IRAN.} \and R. Bakhshandeh
Chamazkoti\thanks{e-mail: ro\_bakhshandeh@mathdep.iust.ac.ir}}
\date{}
\maketitle
\renewcommand{\sectionmark}[1]{}
\begin{abstract} The special linear representation of a compact Lie
group $G$ is a kind of  linear representation of
 compact Lie group $G$ with special properties. It is possible to
 define  the integral of linear
 representation and  extend this concept to special linear
 representation  for next using.
\end{abstract} M.S.C. 2007:  22E46.\\ \hspace{-5mm}{\bf Keywords.}
semi-simple, special linear representation, equivalent, character.
\section{Introduction}
After explicit expressing of special linear representations of
compact Lie group $G$, we consider some of properties  this
structure for compact Lie groups. The first step, in investigating
of compact Lie group, is to obtain a necessary and sufficient
condition in order that a given real Lie algebra $\g$ should be
the Lie algebra of a compact Lie group. It turns out that the
desired criterion consists in the existence on the vector space of
the Lie algebra of a scalar product $\langle u,v\rangle$ which
turns it into a Euclidean space and is invariant with respect to
the adjoint group. Invariance with respect to the adjoint group
amounts to the condition
\begin{equation}\label{1}
\langle[x,u],v\rangle+\langle u,[x,v]\rangle=0~~~;~~~x,u,v\in\g,
\end{equation}
such a Lie algebra will be called a compact Lie algebra. A basic
tool in the entire investigation is the introduction, on an
arbitrary Lie algebra $\g$, of a uniquely determinant scalar
product $\langle u,v\rangle$. If $u,v\in\g$, then scalar product
$\langle u,v\rangle$ is defined to the trace of the linear
transformation $-p_up_v$,
\begin{equation}\label{2}
\langle u,v\rangle=-tr(p_up_v),
\end{equation}
where $p_u$ denotes the inner derivation $p_u=[u,x]$. Clearly
scalar product $\langle u,v\rangle$, according to this definition,
is a symmetric bilinear form on $\g$, also in terms of a
coordinate system, we have
$$\langle u,v\rangle=-{c_\alpha}_j^i{c_\beta}_j^iu^{\alpha}v^{\beta}.$$
An elementary argument shows that a real Lie algebra $\g$ is
compact semi-simple when and only when the scalar product (2)
turns $\g$ into a Euclidean space, i.e., is positive definite.
Since this scalar product is automatically invariant with respect
to all  automorphisms of $\g$, it follows that the group $G_A$ of
all automorphisms of a compact semi-simple Lie algebra $\g$ is a
compact Lie group.  Let $\g$ be a real Lie algebra and let $L$
denote the adjoint group of $\g$. A symmetric bilinear form
$\psi(u,v);u,v\in\g$, defined on $\g$ will be said to be invariant
with respect to $L$ or, simply, invariant if
$$\psi(l(u),l(v))=\psi(u,v),$$
holds for every $l\in L$. Similarly, the associated quadratic form
$\psi(u,v)$ is invariant with respect to the adjoint or, simply,
invariant if
$$\psi(l(u),l(u))=\psi(u,u),$$
holds for arbitrary $l\in L$. Every compact Lie algebra admits a
positive definite quadratic form invariant with respect to the
adjoint group.

Let $\g$ be any real Lie algebra admitting an invariant positive
definite quadratic form. Then the center of $\g$ is precisely the
collection ${\g}_0$ of all elements of $\g$ that are orthogonal to
$\g$ in the sense of the scalar product defined in (2), moreover,
if ${\g}_0$ is trivial then the quadratic form $\psi(u,u)$ of the
scalar product is itself an invariant positive definite quadratic
form. In this case, Lie group $G_A$ of automorphisms of $\g$ is
compact and its Lie algebra $G_A$ coincides with the adjoint
algebra $F$ of $\g$. Finally, if $G$ denotes the component of
identity in $G_A$, then neither $G$ nor $G_A$ contains any central
element other than identity element.
\paragraph{Theorem 1.1.} A
Lie algebra is compact if and only if it admits a positive
definite quadratic form which is invariant with respect to the
adjoint group. (For a proof  and more  details see \cite{[9]}.)

Let $\g$ be a compact semi-simple Lie algebra, and let $G$ be the
component of the identity in the group $G_A$ of all automorphisms
of $\g$, then $G$ is a compact Lie group with trivial center whose
Lie algebra is isomorphic with given algebra $\g$. If $G^{\prime}$
is any Lie group whose Lie algebra is isomorphic with $\g$, then
$G^{\prime}$ is locally isomorphic with $G$; moreover, if
$G^{\prime}$ is a connected global group having trivial center
then it is isomorphic with $G$.
\paragraph{Example 1.2.} Let $A\in
SU(2)$ and ${\goth su}(2)$ be the Lie algebra of $A\in SU(2)$. We
define the adjoint action of $SU(2)$ on ${\goth su}(2)$ by
$$Ad_A(U)=AUA^{-1}=AUA^{*},$$ then each $Ad_A$ is an ${\bf
R}$-linear isomorphism ${\goth su}(2)\longrightarrow {\goth
su}(2)$. We define the scalar product on ${\goth su}(2)$ by
\begin{equation}\label{3}
    \langle X,Y\rangle=-tr(XY)\hspace{1cm}(X,Y\in {\goth su}(2)).
\end{equation}
The elements
\begin{displaymath}
{\hat H} = {1\over\sqrt2}\left( \begin{array}{ccc}
i & 0  \\
0 & i  \\
\end{array} \right),\hspace{1cm}
{\hat E} = {1\over\sqrt2}\left( \begin{array}{ccc}
0 & 1 \\
-1 & 0 \\
\end{array} \right),\hspace{1cm}
{\hat F} = {1\over\sqrt2}\left( \begin{array}{ccc}
0 & i  \\
-i & 0  \\
\end{array} \right)~,
\end{displaymath}
form an orthonormal basis  $\{{\hat H},{\hat E},{\hat F}\}$ of
 ${\goth su}(2)$ with respect to the inner product $\langle~,~\rangle$ ,i.e.,
$$\langle \hat H,\hat H\rangle=\langle\hat E,\hat E\rangle=\langle\hat F,\hat F\rangle=1,$$
$$\langle\hat H,\hat E\rangle=\langle\hat H,\hat F\rangle=\langle\hat E,\hat F\rangle=0.$$
We can define an ${\Bbb R}$-linear isomorphism
$$\theta:{\Bbb R}^3\longrightarrow {\goth su}(2)~~~;~~~\theta(xe_1+ye_2+ze_3)=x{\hat E}+y{\hat H}+z{\hat F},$$
which is also an isometry, i.e.,\\
$$\langle\theta(x),\theta(y)\rangle=x.y\hspace{1cm}(x,y\in{\Bbb R}^3).$$
In particular, the scalar product $\langle~,~\rangle$ defined in
(3),  has positive definite quadratic form on ${\goth su}(2)$,
whence it follows that ${\goth su}(2)$ is semi-simple and compact.
Let $G_A$ denote the group of all automorphisms of ${\goth
su}(2)$. Since an automorphism preserves the scalar product, it
follows that $G_A$ consists  of orthogonal transformations on the
vector space ${\goth su}(2)$ exclusively. Moreover, it is clear
that the orthogonal transformations with positive determinant
preserve the vector product while those with negative determinant
do not. Thus $G_A$ is just the group of all rotations of ${\goth
su}(2)$. In fact, we obtain that Lie algebra ${{\goth su}(2)}_A$
of the group $G_A$ is isomorphic with ${\goth su}(2)$.
\section{Representation theory }
In this section, we introduce several definitions of elementary
representation theory for using  in next
section:\paragraph{Definition 2.1.} A {\it representation} of a
Lie group is defined by a group homomorphism $\rho$ of Lie group
$G$ to $GL(V)$, from Lie group $G$ to the space of invertible
linear transformations on a vector space $V$. If $\rho$ has finite
order, then to each element $x\in G$ there corresponds a matrix
$\rho(x)$, the entries of which we denote by $\rho_j^i(x)$,
$\rho(x)=[\rho_j^i(x)]$, each of the matrix $\rho(x)$ there
corresponds a non-singular linear automorphism of $V$ which we
denote by $\rho_x$. The order of $\rho(x)$ is known as the degree
of $\rho$.

 Let $G$ be a Lie group and
$g\in G$ be an arbitrary element of $G$. A representation
$\rho:G\longrightarrow GL(V)$ is {\it continuous} in $g$ when and
only when for each $\epsilon>0$, there exists a neighborhood $N$
of the identity such that for all $x\in G$, $xg^{-1}\in N$ we have
$||\rho(x)-\rho(g)||<\epsilon$, where $||\rho(x)||$ is the same
suprimum norm on invertible linear transformation $\rho_x$.
\paragraph{Definition 2.2.} Let $\rho:G\longrightarrow GL(V)$ be a
linear representation on a compact Lie group,
$\rho(x)=[\rho_j^i(x)]$. We will say the linear representation has
an {\it invariant integral} when and only when for every
continuous complex function $\rho_j^i$ defined on $G$, there
associated a complex number, denoted by $\int\rho_j^i(x)dx$ and
called integral of $\rho_j^i$ over $G$ such that the following
conditions are satisfied:

(i)~ If $\alpha$ is a complex number then
$$\int\alpha\rho_j^i(x)dx=\alpha\int\rho_j^i(x)dx;$$

(ii)~If $\rho$ and $\xi$ are two linear representations of compact
Lie group $G$ which $\rho(x)=[\rho_j^i(x)]$ and
$\xi(x)=[\xi_j^i(x)]$, and $\rho_j^i$ and $\xi_j^i$ are two
continuous functions then
$$\int(\rho_j^i(x)+\xi_j^i(x))dx=\int\rho_j^i(x)dx+\int\xi_j^i(x)dx;$$

(iii)~ If the function $\rho_j^i$ is everywhere non-negative then
$$\int\rho_j^i(x)dx\geq0;$$
more than above conditions if $\rho_j^i$ everywhere no identically
zero then $\int\rho_j^i(x)dx>0$,

(iv)~ If $\rho_j^i(x)=1$ for every $x\in G$ then
$$\int\rho_j^i(x)dx=1;$$

(v)~ If $a\in G$ be an arbitrary element with $\rho_j^i(a)=1$ then
$$\int\rho_j^i(ax)dx=\int\rho_j^i(xa)dx=\int\rho_j^i(x)dx;$$

(vi)~ $$\int\rho_j^i(x^{-1})dx=\int\rho_j^i(x)dx.$$ Indeed, the
integral of a linear representation is obtained by integrating the
entries of $[\rho_j^i(x)]$ matrix separably, then a matrix with
the same number of rows and columns will be formed.
\paragraph{Definition 2.3.} Two linear representations
$\rho$ and $\xi$ of the same Lie group $G$ are said to be {\it
equivalent} if there exists a constant (independent of $x$) matrix
$A$ such that $\xi(x)=A\rho(x)A^{-1}$ for every $x\in G$.

 Let $U$ be a complex vector
space. A complex valued function $\varphi(x,y)$ of two vectors
$x,y\in U$ is called a {\it Hermitian bilinear form} if

(i)~ $\varphi(\lambda x+\mu
y,z)=\lambda\varphi(x,z)+\mu\varphi(y,z)$, where $\lambda$ , $\mu$
denote complex numbers.

(ii)~ $\varphi(x,y)=\varphi(y,x)$.
 A Hermitian form
$\varphi(x,y)$ is said to be positive definite if $\varphi(x,x)>0$
for $x\neq0$. Clearly any positive definite Hermitian bilinear
form may be taken as scalar product in the space $U$ becomes a
unitary space.

 If $\rho:G\longrightarrow
GL(V)$ is any representation, and $W\subset V$ is an invariant
subspace, meaning that $\rho(g)W\subset W$ for all $g\in G$, then
representation $\rho$ reduces to a subrepresentation of $G$ on
$W$. Trivial subrepresentations are when $W=\{0\}$ or $W=V$. A
{\it reducible} representation is one that contains a nontrivial
subrepresentation; an {\it irreducible} representation, then, is
one that has no nontrivial invariant subspace (subrepresentation).
\paragraph{Example 2.4.} The direct sum of two
representations is reducible since each appears as a
subrepresentation therein. Let $G$ is a compact Lie group, or the
general tensor representations of $GL(n,{\Bbb R})$, all
finite-dimensional reducible representations could be decomposed
into a direct sum of irreducible subrepresentations. However, this
is not uniformly true; a simple counterexample is provided by the
reducible two dimensional representation
\begin{displaymath}
\rho(A) = \left[ \begin{array}{ccc}
1 & \log|\det A|  \\
0 & 1 \\
\end{array} \right],
\end{displaymath}
of the general linear group $GL(n,{\Bbb R})$.\paragraph{Remark I.}
Let $\Delta$ be a reducible set of unitary matrix of order $r$.
Then there exists a single unitary matrix $P$ of order $r$ such
that for every $D\in \Delta$ the matrix ${\tilde D}=PDP^{-1}$ has
the special form
\begin{equation}\label{4}
{\tilde D} = \left[ \begin{array}{ccc}
A & 0  \\
0 & B \\
\end{array} \right]
\end{equation}
where $A$ and $B$ are unitary matrices. We regard the matrices of
$\Delta$ as the matrices of a set of unitary transformations of an
$r$-dimensional unitary space $V$ relative to a fixed orthonormal
basis.
 Let $V$ be an unitary space of finite dimension $m$ and
$G$ be a compact Lie group . A linear representation
$\rho:G\longrightarrow GL(V)$ of $V$ into itself is said to be
{\it unitary} if it preserves the inner product, i.e., for all
$x,y\in V$; $\langle \rho(x),\rho(y)\rangle=\langle x,y\rangle$.
\paragraph{lemma 2.5.} ($Schur^{,}s ~lemma$)Let $\Sigma$ be an irreducible
collection of t -rowed square matrices, let $\Omega$ be an
irreducible set of $n$-rowed square matrices, and $A$ be a
rectangular matrix having $m$ rows and $n$ columns. Suppose that
$\Sigma A=A\Sigma$ then either $A=0$, or else $m=n$ and $A$ is
non-singular.
\paragraph{lemma 2.6.} Let $\Omega$ be an irreducible
set of complex $m\times m$ matrix and $A$ be a square matrix of
order $r$ that commutes with all elements of $\Omega$. Then $A$
has the form $\beta E$ where $\beta$ is a complex number and $E$
is the identity matrix.

\medskip \noindent {\it Proof.} Consider the matrix $B=A-\beta E$
where $\beta$ a complex number chosen, so as to make the
determinant of a zero. Since the $\det(A-\beta E)$ is a polynomial
in $\beta$ with complex coefficients, the existence of such a
number is assured. Moreover, since $A$ commute with all the
elements of $\Omega$ the same true of $B$. Thus we have $\Omega
B=B\Omega$ and according to $Schur^{,}s ~lemma$ (lemma 2.5.), a
must be the zero matrix since $\det(B)=0$ by construction, then
$A=\beta E$.\hfill\ $\triangle$
\section{Special linear representation}
Now we are going  to create a new  concept in linear
representation while perusing compact Lie groups in new
structure.
\paragraph{Definition 3.1.} A complex linear
representation of a compact Lie group $G$, $\Phi:G\longrightarrow
GL(V)$ is called a special linear representation when $\Phi$ has
one and only one distinct equivalent unitary representation.
\paragraph{Note.} Each two special linear
representations that are not-equivalent are called distinct.
\paragraph{Theorem 3.2.} Let $G$ be a compact Lie group
and $V$ be an unitary space. A complex linear representation
$\Phi:G\longrightarrow GL(V)$ is a special linear representation
if only $\Phi_x , \Phi_y ; x,y\in G$, satisfied in positive
definite Hermitian forms for each basis of $V$.

\medskip \noindent {\it Proof.} First we show that
$\Phi$ has an equivalent unitary representation. Let $V$ be an
$m$-dimensional unitary space, where $m$ is degree of $\Phi$, fix
a basis in $V$, and consider the positive definite Hermitian form
\begin{equation}\label{5}
    \psi(u,v)=\sum_{i=1}^mu_i{\overline{v_i}},
\end{equation}
where $u=(u_1,\ldots,u_m)$ , $v=(v_1,\ldots,v_m)$. Substituting
$\Phi_x(u)$ , $\Phi_y(v)$ in (5) we obtain a function
\begin{equation}\label{6}
    \psi_x(u,v)=\psi(\Phi_x(u),\Phi_y(v)),
\end{equation}
which is again a positive definite Hermitian form. Now we define
the new Hermitian form\begin{equation}\label{7}
    \varphi(u,v)=\int\psi_x(u,v)dx.
\end{equation}
This form is again a positive definite. Moreover, taking account
the relation $\Phi_x\Phi_y=\Phi_{xy}$ and the invariance of the
integral, we have
$$\varphi(\Phi_x(u),\Phi_y(v))=\int\psi(\Phi_{xy}(u),\Phi_{xy}(v))dx$$
$$\hspace{11mm}=\int\psi_{xy}(u,v)dx$$
$$\hspace{10mm}=\int\psi_x(u,v)dx$$
$$\hspace{1mm}=\varphi(u,v).$$
This shows that $\varphi(u,v)$ value is invariant under the
substitution of $\Phi_y(u),\Phi_y(v)$ for each $u,v$. We can use
$\varphi(u,v)$ as the inner product in $V$, and select in $V$ a
basis which is orthonormal with respect to $\varphi(u,v)$.
Relative to this basis, $\Phi_y$ corresponds to some matrix
$\tilde\Phi(y)$ and since $\Phi_y$ preserves the inner product it
follows that $\tilde\Phi(y)$ is unitary. Thus $\tilde\Phi$ is
unitary representation of $G$. Denoting by $A$, the matrix
connecting the old basis with the new, we have
$\tilde\Phi(x)=A\Phi(x)A^{-1}$, and therefore for any linear
representation there exists an equivalent unitary representation,
also since only $\Phi_x , \Phi_y$ for $x,y\in G$ satisfy in
positive definite Hermitian forms for each basis of $V$ then,
clearly, $\tilde\Phi(x)$ is one and only one distinct equivalent
unitary representation of $\Phi(x)$ and according to definition
(3.1) $\Phi$ is a special linear representation.\hfill\
$\triangle$
\paragraph{Definition 3.3.} The {\it character}
$\chi(x)$ of a special linear representation $\Phi$ is the trace
of the matrix $\Phi(x)$, $\chi(x)=tr[\Phi_j^i(x)]$. Thus the
character of a special linear representation is a numerical
function defined on $G$. Since $\Phi(x)$ and $A\Phi(x)A^{-1}$ have
the same trace, two equivalent representations have the same
character. The character has the property of invariance, viz,
$$\chi(a^{-1}xa)=tr(\Phi(a^{-1}xa))$$
$$\hspace{27mm}=tr(\Phi(a^{-1})\Phi(x)\Phi(a))$$
$$\hspace{9mm}=tr(\Phi(x))$$
$$\hspace{3mm}=\chi(x).$$
\paragraph{Remark II.} Let $\Phi$ be a reducible  special linear representation of $G$.
According to theorem (3.2) and the remark I, there exists a
matrix $P$ such that the matrices $\xi(x)=P\Phi(x)P^{-1}$ have the
special form
\begin{displaymath}
\xi(x) = \left[ \begin{array}{ccc}
\Phi_1(x) & 0  \\
0 & \Phi_2(x) \\
\end{array} \right],
\end{displaymath}
where $\Phi_1(x)$ and $\Phi_2(x)$ are unitary matrices. We shall
say that splits into the two special representations $\Phi_1$ and
$\Phi_2$. If $\Phi_1$ or $\Phi_2$ should be, themselves, reducible
then they may be split again in the same manner. Thus every
special representations $\xi$ splits into a finite system of
irreducible special representation $\Phi_1 ,\ldots,\Phi_n$. If
$\chi$ denotes the character of $\Phi$ and $\chi_i$ be  the
character of $\Phi_i$ then
$$\chi=\chi_1+\ldots+\chi_n ~.$$
\paragraph{Theorem 3.4.} Let $\Phi$ and $\Psi$ be any
two irreducible distinct special linear representations of $G$,
$\Phi=[\Phi_j^i(x)]$, $\Psi=[\Psi_j^i(x)]$, and let $\chi ,
\chi^{\prime}$ be the characters of $\Phi$ and $\Psi$
respectively. Then the following relation holds:
\begin{equation}\label{}
    \int\chi(x)\chi^{\prime}(x)dx=0.
\end{equation}
\medskip \noindent {\it Proof.} Let $m$ and $n$ be the degrees of the distinct special linear
representations $\Phi$ and $\Psi$, and let $A$ be any constant
$m\times n$ matrix. We define
$$T(x)=\Phi(x)A\Psi(x^{-1}),$$
and
$$T=\int T(x)dx~.$$
We show first that $\Phi(y)T\Psi(x)=T.$ Indeed,
$$\Phi(y)T\Psi(y)=\int \Phi(y)\Phi(x)A\Psi(x^{-1})\Psi(y^{-1})dx$$
$$\hspace{6mm}=\int\Phi(yx)A\Psi((yx)^{-1})dx$$
$$\hspace{-25mm}=T.$$
Thus we have
$$\Phi(x)T=T\Psi(x),$$
for every $x$. By lemma $Schur^{,}s ~lemma$ (2.5) there are only
two possible cases. But if $m=n$ and $T$ is non-singular then we
have
$$\Psi(x)=T^{-1}\Phi(x)T,$$
i.e., the special representations $\Phi$ and $\Psi$ are
equivalent, contrary to hypothesis. Thus $T=0$, it means
\begin{equation}\label{9}
   \int\Phi(x)A\Psi(x^{-1})dx=T=0.
\end{equation}
Now choose the special matrix for $A$ with all entries zero except
for $T$ single one in the $(j,1)$-th entry. Then since
$$\Psi(x^{-1})=\Psi(x)^{-1}=\Psi(x)^{*},$$
relation (9) assumes the form
\begin{equation}\label{10}
   \int\Phi_j^i(x)\overline{\Psi_1^k(x)}dx=0.
\end{equation}
Now since
$$\chi(x)=\sum_{i=1}^m\Phi_i^i(x)~~~,~~~\chi^{\prime}(x)=\sum_{j=1}^n\Psi_j^i(x),$$
then we have relation (8).\hfill\ $\triangle$
\paragraph{Theorem 3.5.} Let $\Phi$ be an irreducible special linear
representation of $G$ of degree $r$, $\Phi(x)=[\Phi_j^i(x)]$, and
$\chi$ be the character of $\Phi$. Then
\begin{equation}\label{11}
 \int\chi(x)\overline{\chi(x)}dx=1.
\end{equation}
\medskip \noindent {\it Proof.} Let $A=[A_j^i]$ be any constant $m\times m$, and define
$$T(x)=\Phi(x)A\Phi(x^{-1}),$$
and
$$T=\int T(x)dx.$$
According the proof of above theorem (3.4) we have
$$\Phi(y)T\Phi(y^{-1})=\int\Phi(y)\Phi(x)A\Phi(x^{-1})\Phi(y^{-1})dx$$
$$\hspace{6mm}=\int\Phi(yx)A\Phi((yx)^{-1})$$
$$\hspace{-22mm}=T.$$
so that $\Phi(x)T=T\Phi(x)$ for arbitrary $x$. Employing remark
(II), we conclude that $T$ has the form $\alpha E$ where $E$ is
the identity matrix while $\alpha$ is a complex number depending
only  on $A$. Thus
\begin{equation}\label{12}
 \int\Phi(x)A\Phi(x^{-1})dx=\alpha E.
\end{equation}
It remains to determine $\alpha$. Then we compute trace on both
side of (12), we have
$$tr\bigg(\int\Phi(x)A\Phi(x^{-1})dx\bigg)=\int tr(\Phi(x)A\Phi(x^{-1}))dx$$
$$\hspace{2cm}=\int tr(A)dx$$
$$\hspace{13mm}= tr(A),$$
and
$$tr(\alpha E)=\alpha r,$$
then
$$\alpha={1\over r}tr(A).$$
In this time, choosing for $A$ the special matrix all of whose
entries are zero except for a single one in the $(j,1)$-th entry,
so that $tr(A)=\delta_j^i$ and
$$\Phi(x^{-1})=\Phi(x)^{-1}=\Phi(x)^{*},$$
we obtain from (12),
\begin{equation}\label{13}
  \int\Phi_j^i(x)\overline{\Phi_1^k(x)}dx={1\over
  r}\delta_k^i\delta_1^j,
\end{equation}
then
\begin{equation}\label{14}
 \int\Phi_j^i(x)\overline{\Phi_j^i(x)}dx={1\over
  r},
\end{equation}
and if $i\neq k$ or $j\neq 1$ then
\begin{equation}\label{15}
 \int\Phi_j^i(x)\overline{\Phi_1^k(x)}dx=0,
\end{equation}
that (14) and (15)  conclude (11) relation immediately .\hfill\
$\triangle$

\end{document}